\documentclass[a4paper]{amsart}
\usepackage[T1]{fontenc}
\usepackage[latin1]{inputenc}
\usepackage{hyperref}
\usepackage{amscd}
\usepackage{epsfig}
\makeatletter
\theoremstyle{plain} 
\newtheorem{thm}{Theorem}[section]
\numberwithin{equation}{section} 
\numberwithin{figure}{section} 
\theoremstyle{plain} 
\newtheorem{cor}[thm]{Corollary} 
\newtheorem{claim}[thm]{Claim} 
\theoremstyle{plain} 
\newtheorem{prop}[thm]{Proposition} 
\theoremstyle{remark}

\theoremstyle{remark}

\usepackage[arrow,matrix,curve,ps]{xy}
\CompileMatrices

\newcommand{\C}{\mbox{$\mathbb{C}$}}

\renewcommand{\O}{\mbox{$\mathcal{O}$}}
\renewcommand{\P}{\mbox{$\mathbb{P}$}}

\DeclareMathOperator{\birat}{bir}

\DeclareMathOperator{\Ext}{Ext}

\DeclareMathOperator{\Hom}{Hom}

\DeclareMathOperator{\Pic}{Pic}

\DeclareMathOperator{\RatCurves}{RatCurves}

\DeclareMathOperator{\Sing}{Sing}

\DeclareMathOperator{\Univ}{Univ}

\makeatother

\begin{document}

\title[Characterizing the projective space]{Characterizing the
 projective space after Cho, Miyaoka and Shepherd-Barron}

\date{May 9, 2001 }

\author{Stefan Kebekus}



\address{Stefan Kebekus, Math.~Institut, Universität Bayreuth, 95440
  Bayreuth, Germany}

\email{stefan.kebekus@uni-bayreuth.de}

\thanks{The author gratefully acknowledges support by the
 Forschungsschwerpunkt ``Globale Methoden in der komplexen Analysis''
 of the Deutsche Forschungsgemeinschaft.}

\urladdr{http://btm8x5.mat.uni-bayreuth.de/\~{}kebekus}

\maketitle

\section{Introduction}

The aim of this paper is to give a short proof of the following
characterization of the projective space.

\begin{thm}\label{thm:main_char}
  Let $X$ be a projective manifold of dimension $n \geq 3$, defined
  over the field $\C$ of complex numbers. Assume that for every curve
  $C\subset X$, we have $-K_X.C \geq n+1$. Then $X$ is isomorphic to
  the projective space.
\end{thm}

A proof was first given in the preprint \cite{CMS00} by K.~Cho,
Y.~Miyaoka and N.~Shepherd-Barron. While our proof here is shorter,
involves substantial technical simplifications and is perhaps more
transparent, the essential ideas are taken from that preprint ---see
section~\ref{sec:attribut}.

This paper aims at simplicity, not at completeness. The methods also
yield other, more involved characterization results which we do not
discuss here. The preprint \cite{CMS00} discusses these thoroughly.

\subsubsection*{Acknowledgement}

This paper was worked out while the author enjoyed the hospitality of
the University of Washington at Seattle, the University of British
Columbia at Vancouver and Princeton University. The author would like
to thank K.~Behrend, J.~Kollár and S.~Kovács for the invitations and
for numerous discussions.

\section{Setup}

\subsection{The space of rational curves}
 
For the benefit of the reader who is not entirely familiar with the
deformation theory of rational curves on projective manifolds, we will
briefly recall the basic facts about the parameter space of rational
curves. Our basic reference is Kollárs book \cite{K96} on rational
curves. The reader might also wish to consider the less technical
overview in \cite{Kebekus-Habil}.

If a point $x \in X$ is given, then there exists a scheme
$\Hom_{\birat}(\P_1,X,[0:1]\mapsto x)$ whose geometric points
correspond to generically injective morphisms from $\P_1$ to $X$ which
map the point $[0:1]\in \P_1$ to $x$. There exists a universal
morphism 
$$
\begin{array}{rccc}
\mu : & \Hom_{\birat}(\P_1,X,[0:1]\mapsto x) \times \P_1 &\to & X \\
& (f,p) & \mapsto & f(p)
\end{array}
$$
If $\mathbb B$ denotes the group of automorphisms of $\P_1$ which
fix the point $[0:1]$, then $\mathbb B$ acts by composition naturally
on the space $\Hom_{\birat}(\P_1,X,[0:1]\mapsto x)$. The quotient in
the sense of Mumford exists. We obtain a diagram as follows.
\begin{equation} \label{diag:rat_curves_x}
  \xymatrix{ 
    \Hom^n_{\birat} (\P_1, X,[0:1]\mapsto x)\times \P_1 \ar[d] 
    \ar[r]^(.65){U_x} \ar@/^.6cm/[rr]^{\mu} &
    {\Univ^{rc}(x,X)} \ar[r]_(.6){ \iota_x} \ar[d]_{ \pi_x} & X \\ 
    \Hom^n_{\birat} (\P_1, X,[0:1]\mapsto x) \ar[r]^(.55){u_x} & 
    {\RatCurves^n(x,X)} }
\end{equation} 
Here $\Hom^n_{\birat}(\ldots)$ is the normalization of
$\Hom_{\birat}(\ldots)$, the morphisms $U_x$ and $u_x$ have the
structure of principal $\mathbb B$-bundles and $\pi_x$ is a
$\P_1$-bundle. The restriction of $\iota_x$ to any fiber of $\pi_x$ is
generically injective, i.e.~birational onto its image.

The space $\RatCurves^n(x,X)$ is called the ``space of rational curves
through $x$''. This name is perhaps a bit misleading because the
correspondence
$$
\begin{array}{rccc}
  e: & \RatCurves^n(x,X) & \to & \{ \text{rational curves in $X$ which
    contain $x$}\} \\
  & h & \mapsto & \iota_x( \pi^{-1}(h))
\end{array}
$$
is not bijective in general. Although $e$ is surjective, it may
happen that the restriction of $e$ to an irreducible component $H
\subset \RatCurves^n(x,X)$ is only generically injective: several
points in $H$ may correspond to the same rational curve.

\subsection{Results of Mori's theory of rational curves}

The following theorem summarizes some of the classic results of Mori
theory, in particular Mori's famous existence theorem for rational
curves on manifolds where $K_X$ is not nef. While most statements can
been found explicitly or implicitly in the papers \cite{Mor79} and
\cite{KMM92}, some results found their final formulation only years
later. We refer to \cite{K96} for proper attributions.

\begin{thm}[Classic results on families of rational curves]\label{thm:Mori-thy}
  Under the assumptions of theorem~\ref{thm:main_char}, let $x$ be a
  very general point\footnote{A ``very general point'' is, by
    definition, a point which is not contained in a countable union of
    certain hypersurfaces.} of $X$. Then \cite[thm.~II.5.14]{K96}
  there exists a rational curve $\ell \subset X$ which contains $x$
  and satisfies $-K_X.\ell = n+1$.
 
  Let $ H_x \subset \RatCurves^n(x,X)$ be the irreducible component
  which contains the point corresponding to $\ell$ and consider the
  restriction of the diagram~(\ref{diag:rat_curves_x}) above:
  \begin{equation}
    \label{eq:eval_and_proj}
    \xymatrix{ { U_x} \ar[r]^{ \iota_x}
      \ar[d]^{ \pi_x}_{\txt{\scriptsize $\P_1$-bundle}} & {X} \\
      { H_x} }
  \end{equation}
  Then the following holds.
  \begin{enumerate}
  \item The variety $ H_x$ is compact \cite[prop.~II.2.14]{K96},
    smooth \cite[cor.~II.3.11.5]{K96} and has dimension $\dim H_x =
    n-1$ \cite[thms.~II.3.11 and II.1.7]{K96}.
  \item The evaluation morphism $ \iota_x$ is finite away from
    $\iota_x^{-1}(x)$ (Mori's Bend-and-Break,
    \cite[cor.~II.5.5]{K96}). In particular, $\iota_x$ is surjective.
  \item If $\ell \subset X$ is curve corresponding to a general point
    of $ H_x$, then $\ell$ is smooth \cite[thm.~II.3.14]{K96} and the
    restriction $T_X|_\ell$ is an ample vector bundle on $\ell$
    \cite[cor.~II.3.10.1]{K96}.
 \end{enumerate}
\end{thm}

\subsection{Singular rational curves}

It was realized very early by Miyaoka (\cite{Miy92}, see also
\cite[prop.~V.3.7.5]{K96}) that the singular curves in the family
$H_x$ play a pivotal rôle in the characterization problem. A
systematic study of families of singular curves, however, was not
carried out before the paper \cite{Keb00a}. In that paper, the author
gave a sharp bound on the dimension of the subvariety $H^{\Sing}_x
\subset H_x$ whose points correspond to singular rational curves and
described the singularities of those curves which are singular at $x$.

The following theorem summarizes the results of \cite{Keb00a} which
form the centerpiece of our argumentation. A singular curve is called
``immersed'' if the normalization morphism has constant rank one. A
singular curve which is not immersed is often said to be ``cuspidal''.

\begin{thm}[Singular curves in $H_x$, {\cite[thm.~3.3]{Keb00a}}]
  \label{thm:bounds_for_sing}
  The closed subfamily $H^{\Sing}_x \subset H_x$ of singular curves
  has dimension at most one. The subfamily $H^{\Sing,x}_x \subset
  H^{\Sing}_x$ of curves which are singular at $x$ is at most finite.
  If $H^{\Sing,x}_x$ is not empty, then the associated curves are
  immersed.
\end{thm}

In our setup, we obtain a good description of $\iota_x^{-1}(x)$ as an
immediate corollary.

\begin{cor}\label{cor:section_and_points}
  The preimage $ \iota^{-1}_x(x)$ contains a section, which we call
  $\sigma_\infty$, and at most a finite number of further points,
  called $z_i$. The tangential morphism $T\iota_x$ has rank one along
  $\sigma_\infty$.
\end{cor}

The universal property of the blow-up \cite[prop.~II.7.13]{Ha77}
therefore allows us to extend diagram~(\ref{eq:eval_and_proj}) as
follows:
$$
\xymatrix{ & & {\hat X} \ar[d]^{\txt{\scriptsize blow-up of $x$}}
 \\ { U_x} \ar[rr]^{ \iota_x} \ar[d]^{
 \pi_x}_{\txt{\scriptsize $\P_1$-bundle}} \ar@{-->}@/^/[rru]^{\hat
 \iota_x} & & {X}
 \\
 { H_x} }
$$
where the rational map $\hat \iota_x$ is well-defined away from the
points $z_i$. 

We end with a further description of $\hat \iota_x$.

\begin{prop}[{\cite[thm.~3.4]{Keb00a}}]\label{prop:finiteness}
  If $E \cong \P(T_X^*|_x)$ is the exceptional divisor of the
  blow-up\footnote{We use Grothendieck's notation: if $V$ is a vector
    space, then $\P(V^*) = V\setminus \{0\}/\C^*$. This drives the
    aficiónados of the older literature to the drink.}, then the
  restricted morphism
  $$
  \hat \iota_x|_{\sigma_\infty} : \sigma_\infty \to E
  $$
  is finite. In particular, since $\dim \sigma_\infty = \dim H_x =
  n-1$, the morphism $\hat \iota_x|_{\sigma_\infty}$ is surjective.
\end{prop}

\section{Proof of the characterization Theorem}

\subsection{The neighborhood of $\sigma_\infty$}
\label{sec:step1}

As a first step towards the proof of theorem~\ref{thm:main_char}, we
need to study the neighborhood of the section $\sigma_\infty \subset
 U_x$.

\begin{prop}\label{prop:neighborhood}
  If $E \cong \P(T_X^*|_x)$ is the exceptional divisor of the blow-up, then
  \begin{enumerate}
  \item The restricted morphism $\hat \iota_x|_{\sigma_\infty}$ is an
    embedding. In particular, $H_x \cong \sigma_\infty \cong \P_{n-1}$.
  \item The tangent map $T \hat \iota_x$ has maximal rank along
    $\sigma_\infty$. In particular, $N_{\sigma_\infty,U_x} \cong
    N_{E,\hat X} \cong \O_{\P_{n-1}}(-1)$.
  \end{enumerate}
\end{prop}

The remaining part of the present section~\ref{sec:step1} is devoted
to a proof of proposition~\ref{prop:neighborhood}.  Note that
statement~(2) follows immediately from statement~(1) and from
corollary~\ref{cor:section_and_points}. To show statement~(1) requires
some work.

By proposition~\ref{prop:finiteness} and by Zariski's main theorem, we
are done if we show that $\hat \iota_x|_{\sigma_\infty}$ is
birational. Assume for the moment that $\hat \iota_x|_{\sigma_\infty}$
is \emph{not} birational, let $\ell \subset X$ be a general curve
associated with $H_x$ and let $F \subset U_x$ be the corresponding
fiber of $ \pi_x$. The subvariety $ \iota_x^{-1}(\ell) \subset U_x$
will then contain a curve $B$ such that
\begin{enumerate}
\item $B$ is not contained in $\sigma_\infty$.
\item $B$ is not a fiber of the projection $\pi_x$.
\item $B\cap \sigma_\infty$ contains a point $y_1$ which is different
 from $y_0 := F \cap \sigma_\infty$.
\end{enumerate}
In order to see that we can find a curve $B$ which is not a fiber of
the projection $\pi_x$, recall that the correspondence between points
in $H_x$ and curves in $X$ is generically injective and that $\ell$
was generically chosen.

Summing up, in order to show proposition~\ref{prop:neighborhood}, it
suffices to show the following claim.
\begin{claim}\label{claim1}
  Let $\ell \subset X$ be a general curve associated with $H_x$ and
  let $B \subset \iota_x^{-1}(\ell)$ be any curve which satisfies
  items (1) and (2). Then $B$ is disjoint from $\sigma_\infty$.
\end{claim}
A proof will be given in the next few subsections.

\subsubsection{The normal bundle of $\ell$}

Since $\dim H_x = n-1 > 1$, it is a direct consequence of
theorem~\ref{thm:bounds_for_sing} that $\ell$ is smooth and therefore
isomorphic to the projective line. A standard theorem, which is
attributed to Grothendieck, but probably much older, asserts that a
vector bundle on $\P_1$ always decomposes into a sum of line bundles.
For the restriction of the tangent bundle $T_X$ to $\ell$, all
summands must be positive by theorem~\ref{thm:Mori-thy}(3). The
splitting type is therefore known:
$$
T_X|_\ell \cong \O(2)\oplus\O(1)^{\oplus n-1}.
$$
The normal bundle of $\ell$ in $X$ is thus isomorphic to
\begin{equation}
 \label{eq:normal_bdle}
 N_{\ell/X} \cong \O(1)^{\oplus n-1}. 
\end{equation}
We will use this splitting later in section~\ref{sec:self-intersect_nmbrs}
to give an estimate on certain self-intersection numbers.

\subsubsection{Reduction to a ruled surface}

As a next step let $\tilde B$ be the normalization of $B$ and perform
a base change via the natural morphism $\tilde B \to H_x$. We obtain a
diagram as follows:
$$
\xymatrix{ { U_B} \ar@/^0.6cm/[rrr]^{ \iota_B}
  \ar[rr]^{\gamma}_{\txt{\scriptsize finite base change}}
  \ar[d]_{\txt{\scriptsize $\pi_B$\\ \scriptsize $\P_1$-bundle}} & & {
    U_x} \ar[r]^{ \iota_x} \ar[d]^{ \pi_x}
  & {X} \\
  {\tilde B} \ar[rr] & & { H_x} }
$$
The bundle $ U_B$ will now contain two distinct distinguished
sections. Let $\sigma_{B,\infty} \subset \iota_B^{-1}(x)$ be
the section which is contracted to a point and choose a component
$\sigma_{B,0} \subset \gamma^{-1}(B)$. In order to prove
claim~\ref{claim1} we have to show that these sections are disjoint.
\begin{figure}
 \scriptsize
 $$
 \xymatrix{ {
 \begin{picture}(4,4)(0,0)
  \put(0.00, 0.00){\epsfxsize 4cm \epsffile{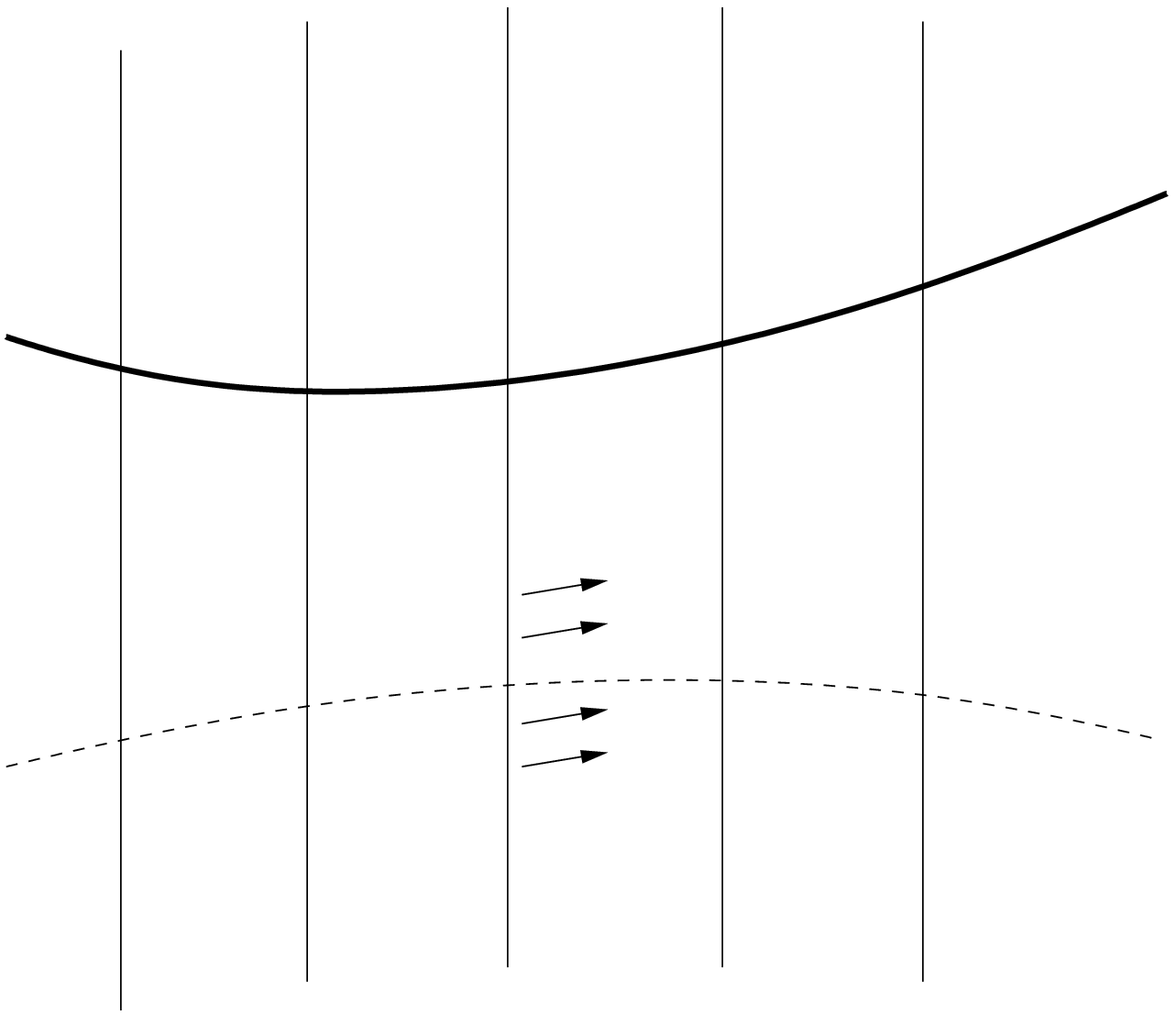}}
  \put(4.10, 2.70){$\sigma_{B,0}$}
  \put(1.75, 3.40){$\ell_y$}
  \put(3.20, 3.30){$\ell_x$}
  \put(4.10, 0.90){$\sigma_{B,\infty}$}
  \put(3.20, 2.30){$z_i$}
  \put(3.07, 2.41){$\bullet$}
  \put(-0.3, 3.7){\normalsize $U_B$}
 \end{picture}} \ar[rr]^(.45){\iota_B} \ar[d]_(.75){\pi_B} & & {
  \begin{picture}(5,4)(-0.2,0)
  \put(0.0, 0.0){\epsfxsize 5cm \epsffile{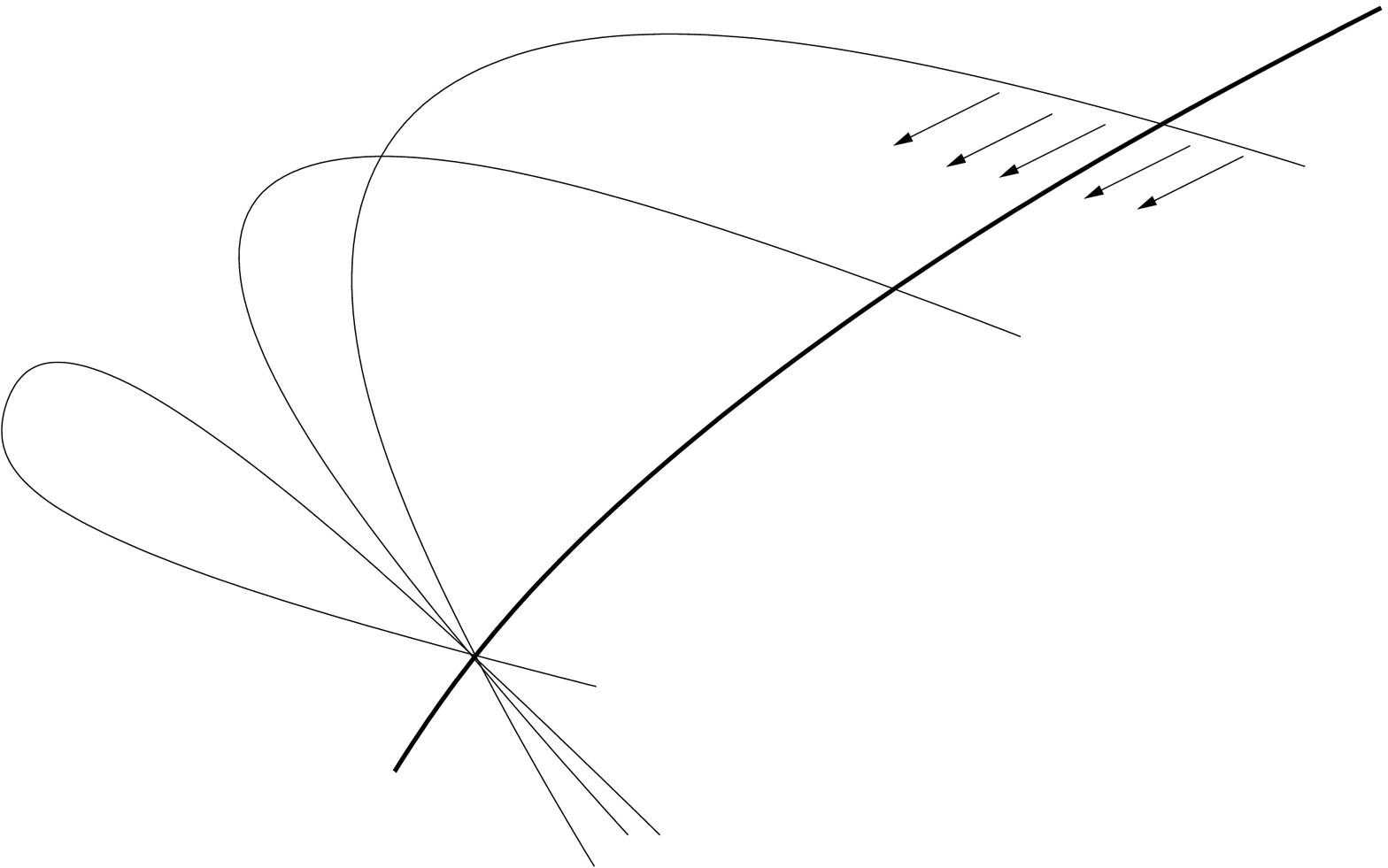}}
  \put(5.0, 3.2){$\ell$}
  \put(0.0, 1.0){$\ell_x$}
  \put(1.8, 3.3){$\ell_y$}
  \put(4.1, 2.9){$y$}
  \put(4.11, 2.62){$\bullet$}
  \put(3.15, 2.02){$\bullet$}
  \put(1.6, 0.4){$x$}
  \put(1.65, 0.7){$\bullet$}
  \put(4.0, 1.0){\normalsize $X$}
  \end{picture}
  } \\ {
  \epsfxsize 4cm \epsffile{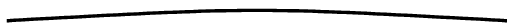}
  \begin{picture}(0,0)(-0.3,0)~
  \put(-5, 0.0){\normalsize $\tilde B$}
  \end{picture}
  } }
 $$
\caption{Reduction to a ruled surface}
\end{figure}

\subsubsection{Estimate for the self-intersection of $\sigma_{B,0}$}
\label{sec:self-intersect_nmbrs}

Let $d$ be the mapping degree of the restricted evaluation
$\iota_B|_{\sigma}$. We will now show that the self-intersection
number of the distinguished section $\sigma_{B,0}$ in the ruled
surface $U_B$ is at most $d$.

\begin{claim}
  The natural map
  $$
  T \iota_B: N_{\sigma_{B,0}, U_B} \to N_{\ell,X}
  $$
  between the normal bundles is not trivial.  
\end{claim}
\begin{proof}
  Let $\hat H_x \subset H_x$ be the closed proper subvariety whose
  points correspond to curves which are either not smooth or whose
  normal bundle is not of the form~(\ref{eq:normal_bdle}). Since
  $\ell$ was generically chosen, $\ell$ is not contained in the proper
  subvariety $\iota_x(\pi_x^{-1}(\hat H_x)) \subset X$. Consequence:
  if $F \subset U_B$ is a general fiber of the morphism $\pi_B$, then
  $\iota_B(F)$ is a smooth curve with normal bundle $N_{\iota_B(F),X}
  \cong \O(1)^{\oplus n-1}$, and the tangent map $T\iota_B$ has rank
  two along $F \setminus (F\cap \sigma_{B,\infty})$. In particular,
  $\iota_B$ has maximal rank at $F\cap \sigma_{B,0}$, and the claim is
  shown.
\end{proof}

We obtain the estimate
\begin{equation}
 \label{eq:estimate_for_B}
 \sigma_{B,0}^2 = \deg N_{\sigma_{B,0}, U_B} \leq d \cdot
 \underbrace{(\text{max. degree of sub-linebundles in $N_{\ell,X}$})}_{=1
   \text{ by equation~(\ref{eq:normal_bdle})}} = d. 
\end{equation}

\subsubsection{Intersection numbers on the ruled surface}
\label{sec:intersect_nmbrs}

Let $F$ be a fiber of the projection $ \pi_B$ and $H\in \Pic X$
be any ample line bundle. We obtain the following list of
intersection numbers.
\begin{align*}
 \iota_B^*(H).\sigma_{B,\infty} & = 0 &&
 \text{because $\sigma_{B,\infty}$ is contracted to a point}\\
 \iota_B^*(H).\sigma_{B,0}  & = d \cdot 
 \iota_B^*(H).F && \text{$ \iota|_F:F\to \ell$ is
 birational and $ \iota|_{\sigma_{B,0}}$ is $d:1$}\\
 \sigma_{B,0}.F & = 1 && \text{because $\sigma_{B,0}$ is a section}
\end{align*}
Consequence: we may write the following numerical equivalence of
divisors on $ U_B$:
$$
\sigma_{B,0} \equiv \sigma_{B,\infty} + d\cdot F.
$$
We end the proof of claim~\ref{claim1} and of
proposition~\ref{prop:neighborhood} with the calculation
\begin{align*}
  \sigma_{B,0}^2 & = \sigma_{B,0}\cdot (\sigma_{B,\infty}+d\cdot F) \\
  & = \sigma_{B,0} \cdot \sigma_{B,\infty} + d.
\end{align*}
The inequality~(\ref{eq:estimate_for_B}) shows that $\sigma_{B,0}
\cdot \sigma_{B,\infty} = 0$. The distinguished sections are therefore
disjoint. The proof of proposition~\ref{prop:neighborhood} is
finished.

\subsection{Factorization of $\iota_x$}

To end the proof of theorem~\ref{thm:main_char}, consider the
Stein-factorization of the morphism $ \iota_x$. We obtain a sequence
of morphisms
$$
\xymatrix{ { U_x} \ar@/^/[rr]^{ \iota_x} 
 \ar[r]_{\alpha} & {Y} \ar[r]_{\beta} & X}
$$
where $\alpha$ contracts the divisor $\sigma_\infty$, and $\beta$
is a finite map. 

Since $R^1{\pi_x}_*(\O_{U_x}) = 0$, the push-forward of the twisted
ideal sheaf sequence
$$
0 \to \O_{U_x} \to \O_{U_x}(\sigma_\infty) \to
\underbrace{\O_{U_x}(\sigma_\infty)|_{\sigma_\infty}}_{\cong
  \O_{\P_{n-1}}(-1)} \to 0
$$
gives a sequence
$$
0 \to \O_{\P_{n-1}} \to \mathcal E \to \O_{\P_{n-1}}(-1) \to 0
$$
on $H_x \cong \P_{n-1}$ where $\mathcal E$ is a vector bundle of
rank two and $U_x \cong \P(\mathcal E^*)$. Since
$\Ext^1_{\P_{n-1}}(\O_{\P_{n-1}}(-1),\O_{\P_{n-1}}) = 0$, the bundle
$U_x$ is thus isomorphic to
$$
U_x \cong \P(\O_{\P_{n-1}}(-1)\oplus \O_{\P_{n-1}}).
$$
Consequence: there exists a morphism $\alpha':U_x \to \P_n$ which
contracts $\sigma_\infty$. An application of Zariski's main theorem
shows that $\alpha = \alpha'$. In particular, $Y \cong \P_n$.

The fact that $\beta$ is an isomorphism now follows from \cite{Laz84}:
note that the $\beta$-images of the lines through $\alpha(x)$ are the
curves associated with $H_x$. Recall the adjunction formula for a
finite, surjective morphism:
$$
-K_{\P_n} = \beta^*(-K_X)+(\text{branch divisor})
$$
To see that $\beta$ is birational, and thus isomorphic, it is
therefore sufficient to realize that for a curve $\ell \in H_x$, we
have 
$$
-K_X.\ell = -K_{\P_n}.(\text{line}) = n+1.
$$
This ends the proof of theorem~\ref{thm:main_char}.

\subsection{Attributions}
\label{sec:attribut}

The reduction to a ruled surface and the calculation of the
intersection numbery have already been used in \cite{Miy92} (see also
\cite[prop.~V.3.7.5]{K96}) to give a criterion for the existence of
singular rational curves. The estimate~(\ref{eq:estimate_for_B}) is
taken from \cite{CMS00} where a similar estimate is used in a more
complex and technically involved situation to prove a statement
similar to proposition~\ref{prop:neighborhood}.

The calculation that $U_x \cong \P(\O(-1)\oplus \O)$ is modelled after
\cite{Mor79}.

\bibliographystyle{alpha}

\end{document}